\documentclass[journal]{IEEEtran}
\usepackage{xcolor,soul,framed} 
\colorlet{shadecolor}{yellow}
\usepackage[pdftex]{graphicx}
\graphicspath{{../pdf/}{../jpeg/}}
\DeclareGraphicsExtensions{.pdf,.jpeg,.png}
\usepackage[cmex10]{amsmath}
\usepackage{array}
\usepackage{mdwmath}
\usepackage{mdwtab}
\usepackage{eqparbox}

\usepackage{lineno}
\usepackage{amsmath}
\usepackage{graphicx}
\usepackage{float}
\usepackage{cite}
\usepackage{amssymb}
\usepackage{amsthm}
\usepackage{subfigure}
\usepackage{epsfig}
\usepackage{multirow}
\usepackage{pdfpages}
\usepackage{float}
\usepackage{nomencl}
\usepackage{subfigure}
\usepackage{epstopdf}
\usepackage{color}
\usepackage{url}
\usepackage{tabularx}
\usepackage{hhline}
\usepackage{rotating}
\usepackage{fancyhdr}
\usepackage[marginal]{footmisc}
\usepackage{microtype}
\usepackage{graphicx}
\usepackage{subfigure}
\usepackage{booktabs} 
\usepackage{amsmath}
\usepackage{multirow}
\usepackage{threeparttable}
\usepackage{stfloats}
\usepackage{makecell}

\usepackage[ruled,linesnumbered]{algorithm2e}

\usepackage{autobreak}

\usepackage{nomencl}

\usepackage{algorithmic}

\allowdisplaybreaks[1] 

\begin{document}
\title{Improved Successive Branch Reduction for Stochastic Distribution Network Reconfiguration}
\author{Wanjun~Huang,~Changhong~Zhao,~\textit{Senior Member},~\textit{IEEE}
\thanks{This work was supported by the CUHK faculty startup project No. 4930935 and the Hong Kong Rsearch Grants Council through ECS Award No. 24210220. (Corresponding author: Changhong Zhao.)}
\thanks{W. Huang and C. Zhao are with the Department of Information Engineering, the Chinese University of Hong Kong, New Territories, Hong Kong SAR, China. Emails: \{wjhuang, chzhao\}@ie.cuhk.edu.hk}
\thanks{This work has been submitted to the IEEE for possible publication. Copyright may be transferred without notice, after which this version may no longer be accessible.}
}
\maketitle

\newtheorem{thm}{Theorem} 
\newtheorem{definition}{Definition}
\newtheorem{lem}[thm]{Lemma}
\newtheorem{assumption}{Assumption}

\begin{abstract}
We propose an improved successive branch reduction (SBR) method to solve stochastic distribution network reconfiguration (SDNR), a mixed-integer program that is known to be computationally challenging. 
First, for a special distribution network with a single redundant branch, we propose an improved design for a one-stage SBR algorithm in the literature to incorporate uncertain renewable generations and loads. Based on solving stochastic optimal power flow, the improved algorithm identifies and searches through a small set of candidate branches, from which it determines the optimal branch to open and obtains a radial network with the minimum expected operational cost. Then, for a general network with multiple redundant branches, we design a heuristic two-stage SBR algorithm based on a close-and-open procedure that iteratively runs the proposed one-stage SBR algorithm. Numerical results on the IEEE 33-bus and 123-bus distribution network models verify the proposed method in terms of optimality and computational efficiency.
\end{abstract}
\begin{IEEEkeywords}
Stochastic distribution network reconfiguration, successive branch reduction, optimal power flow.
\end{IEEEkeywords}

\section*{Nomenclature}
\addcontentsline{toc}{section}{Nomenclature}
\begin{IEEEdescription}[\IEEEusemathlabelsep\IEEEsetlabelwidth{$q_1$, $q_2$,\quad}]
	\item[\textbf{Sets}]
	\item[$\mathcal{N}$] The set of buses in a network, including substation buses $\mathcal{N}_s$ and non-substation buses $\mathcal{N}_d$.
	\item[$\mathcal{E}$] The set of branches $e=ij=(i,j)\in\mathcal{E}$.
	\item[$\mathcal{W}$] The set of scenarios $w\in\mathcal{W}$ for uncertain renewable generations and loads. 
	\item[$\mathcal{A}$] The set of feasible switch status vectors $\boldsymbol{\alpha}\in\mathcal{A}$ that lead to radial networks.
	\item [$\mathcal{P}$] The set of branches that form a loop.
	\item [$\mathcal{K}$] The set of candidate branches to open.
	\item[] 
	\item[\textbf{Variables}] 
	\item[$p_{i}^{w}$, $q_{i}^{w}$] The active and reactive power injections at bus $i$ in scenario $w$.
	\item[$\hat{p}_{ri}^{w}$, $\hat{q}_{ri}^{w}$] The active and reactive renewable power generations at bus $i$ in scenario $w$. 
	\item[$\hat{p}_{di}^{w}$, $\hat{q}_{di}^{w}$] The active and reactive power loads at bus $i$ in scenario $w$.
	\item[$p_{ij}^w$, $q_{ij}^{w}$] The active and reactive power flows on branch $ij$ in scenario $w$.
	\item[$v_{i}^{w}$] The squared voltage magnitude at bus $i$ in scenario $w$.
	\item[$l_{ij}^{w}$] The squared current magnitude on branch $ij$ in scenario $w$.
	\item[$\boldsymbol{x}$] The continuous decision variables:
	\begin{eqnarray}
	    \boldsymbol{x}&=&\left(\boldsymbol{x}^{w},\forall w\in\mathcal{W}\right)\nonumber\\
	&=& \left(p_{i}^{w}, q_{i}^{w}, \forall i\in \mathcal{N}_s;~v_{i}^{w},\forall i\in\mathcal{N};\right.\nonumber\\
	&& \left.p_{ij}^{w},q_{ij}^{w},l_{ij}^{w},\forall ij\in \mathcal{E}; ~\forall w\in\mathcal{W}\right).\nonumber
	\end{eqnarray}
	
	\item[$\widetilde{p}_{ij}$] The expected active power flow on branch $ij$.
	
	\item[$\widetilde{p}_{i}^{\mathcal{P}}$] The expected active power injection from bus $i$ into loop $\mathcal{P}$. 
	
	\item[$\widetilde{p}_0(e)$] The objective value of SDNR attained by opening a single redundant branch $e$.
	
	\item[$\alpha_{ij}$] The binary variable indicating the switch status on branch $ij$, collected in $\boldsymbol{\alpha}=(\alpha_{ij},\forall ij\in\mathcal{E})$.
	\item[$\mathcal{G}(\boldsymbol{\alpha})$] The network under switch status $\boldsymbol{\alpha}$.
	\item[] 
	\item[\textbf{Parameters}] 
	
	\item[$\boldsymbol{\pi}$] The probability distribution $\boldsymbol{\pi}=(\pi_w,\forall w\in\mathcal{W})$ of the uncertainty scenarios. 
	\item[$L$] The number of redundant branches (chordless loops) in a network.
	\item[$n_r$] The number of active power-injecting buses in a loop considered by the proposed algorithm.
	\item[$\textnormal{r}_{ij}$, $\textnormal{x}_{ij}$] The resistance and reactance of branch $ij$.
\end{IEEEdescription}

\section{Introduction}
Distribution network reconfiguration (DNR) is an important technology to improve the energy efficiency of distribution networks without extra equipment investment \cite{baran1989network}. It mainly aims to minimize the network operational cost (e.g., power loss) by changing the open/closed status of the switches on the branches (power lines). Mathematically, DNR is a mixed-integer nonlinear program that is difficult to solve due to the nonconvex power flow constraints and binary switching actions subject to the requirement for a radial (tree) topology \cite{lavorato2011imposing,wang2020radiality}. 

The existing methods for DNR mostly fall into three categories: mathematical programming, heuristics or meta-heuristics, and machine learning \cite{mahdavi2021reconfiguration,yin2020data}. 
Various mathematical programming techniques have been developed to improve the computational efficiency of DNR, including Benders decomposition \cite{Khodr2009distribution}, column-and-constraint generation \cite{zheng2020deep,zheng2021adaptive}, linearization into a mixed-integer linear program \cite{Francisco2012mixed}, and convex relaxation into a mixed-integer quadratic, quadratically constrained, or conic program \cite{Taylor2012convex,Jabr2021minimum,mahdavi2021reconfiguration}. 
However, the computation process of mathematical programming can be time-consuming for large systems \cite{mahdavi2021reconfiguration}. Common tricks to alleviate computational burdens, such as linearization, often sacrifice the accuracy and optimality of solutions. 

As the second category of DNR methods, heuristic algorithms exploit the topological features of a network to solve the DNR problem at lower computational complexity. 
The commonly used algorithms include iterative branch exchange \cite{civanlar1988distribution}, successive branch reduction (SBR) \cite{gomes2006new}, and switch opening and exchange \cite{zhan2020switch}. The advantage of those heuristic algorithms is that their solutions are generally feasible in terms of the physical laws and operational limits of the network. Starting from a radial topology, the iterative branch exchange algorithm opens a closed switch and closes an open switch in each iteration, until the solution cannot be further improved \cite{civanlar1988distribution}. It is sensitive to the initialization and may take many iterations to converge. The SBR algorithm assumes all the switches are initially closed and then opens them sequentially based on a certain rule until a radial topology is obtained \cite{gomes2006new}. It thus only needs a trivial initialization and a number of iterations bounded by the number of redundant branches. Based on convex relaxation of optimal power flow (OPF), a heuristic algorithm was developed in \cite{peng2014feeder} to improve the computational efficiency of SBR. The switch opening and exchange algorithm combines the iterative branch exchange and SBR methods to provide a more accurate solution given that a higher computational overhead is affordable \cite{zhan2020switch}. Compared to the heuristic algorithms above, metaheuristics such as genetic algorithms \cite{jakus2020optimal} and particle swarm optimization \cite{pegado2019radial} are problem-independent and more general. Nevertheless, they suffer heavier computational burdens from their random selection processes and often end up with inconsistent solutions from different runs. Moreover, their performance relies heavily on parameter tuning, and they may have difficulty in satisfying the radiality condition. Reference \cite{mahdavi2021reconfiguration} provided a detailed review of the methods above and their pros and cons.

Different from the categories above, machine learning methods extract optimal operational knowledge from historical data without needing an accurate physical model of the network. There are mainly two types of machine learning methods: supervised learning and unsupervised learning (e.g., reinforcement learning). A supervised learning approach may apply deep neural networks to learn the relationship between the system state and the optimal topology \cite{huang2021distribution,han2022learning}, which requires a large data set. The reinforcement learning approach learns the optimal control policy by directly interacting with the real or a simulated physical environment  \cite{gao2020batch,wang2021distribution}. Its performance is largely dependent on the hyper-parameters selected from experience. In the machine learning methods, the online training process may contain safety risks, while the offline training may not provide accurate solutions. Another major concern about the machine learning methods is that the feasibility of the solution may not be guaranteed.

The basic, deterministic version of DNR assumes generations and loads are fixed. The more complicated robust DNR \cite{Lee2015robust,zheng2020deep} and \textit{stochastic DNR (\textbf{SDNR})} \cite{Dorostkar2019stochastic,zhan2020switch} are formulated to deal with the uncertainty in renewable generations (e.g., solar and wind) and loads. The robust DNR aims to minimize the operational cost in the worst-case scenario of the uncertain quantities. Therefore, the solution of robust DNR may be too conservative to be economically efficient. 
The SDNR considers the statistical properties of the uncertain quantities and minimizes the expected cost over all the possible scenarios, thus being less conservative and more economical than the robust formulation. The main challenge to SDNR lies in its computational complexity under a large number of scenarios.  

To overcome this challenge, we develop an improved SBR method. The proposed method is inspired by a baseline SBR algorithm in \cite{peng2014feeder}, which has low computational complexity, small optimality loss, and no need for parameter tuning and special initialization. However, the algorithm in \cite{peng2014feeder} only applies to deterministic DNR with fixed positive loads. Our method improves the one in \cite{peng2014feeder} in the following aspects:
\begin{itemize}
	\item For a network with a single redundant branch (compared to a radial network), an improved one-stage SBR algorithm is developed to determine the optimal branch to open in an SDNR problem. The proposed algorithm is built on solving the second-order cone relaxation of stochastic OPF problems. Compared to the algorithm in \cite{peng2014feeder} that is applicable to deterministic DNR with positive loads only, the improved algorithm can incorporate uncertain renewable generations and loads while retaining the optimality of the solution under certain conditions.  
	\item For a network with multiple redundant branches, a two-stage SBR algorithm is developed based on the one-stage algorithm above. In the first stage, from an initial condition where all the branches are closed, a set of redundant branches are selected to open, leading to a radial network. In the second
	stage, each of the branches opened in the first stage is iteratively closed to create a network with a single redundant branch. In each iteration, the one-stage algorithm is used to find the optimal branch to open. The objective values are compared over all such iterations to determine the ultimate radial network with the minimum expected operational cost.
\end{itemize}

The proposed two-stage algorithm needs to solve only one stochastic OPF in its first stage and check a small set of candidate branches in each iteration of the second stage. The stochastic OPF can be further decomposed into multiple deterministic OPF problems solved in parallel. Therefore, it is computationally more efficient than the switch opening and exchange method \cite{zhan2020switch} that solves a number of stochastic OPF problems, each corresponding to a different topology. The optimality and computational efficiency of the proposed method is verified via numerical simulations of the IEEE 33-bus and 123-bus distribution networks, under varying penetration levels of renewable generations and numbers of uncertainty scenarios. 

The remainder of the paper is organized as follows. Section II introduces the model and problem formulation of SDNR. Our one-stage and two-stage SBR algorithms are elaborated in Section III. The numerical case studies are presented in Section IV. Section V provides our conclusion.

\section{Model and Problem Formulation}
\subsection{Stochastic Distribution Network Reconfiguration}
Due to the concern for the lifetime of switches, it is impractical to change the switch status as frequently as the fast fluctuations of renewable generations and loads, especially those of the growing solar photovoltaic generations in distribution networks. To address this concern, the stochastic distribution network reconfiguration (SDNR) aims to determine an optimal radial topology of the network to achieve the minimum expected operational cost over the possible scenarios of uncertain generations and loads. 

Consider a distribution network with a set of buses $\mathcal{N}$ and a set of branches $\mathcal{E}$. Each branch is arbitrarily assigned a reference direction, say from bus $i$ to bus $j$, and is represented as $ij\in\mathcal{E}$ or $(i,j)$ interchangeably.
The bus set $\mathcal{N}$ is divided into substation buses $\mathcal{N}_{s}$ and non-substation buses $\mathcal{N}_{d}$.
The substation buses are connected to the main grid outside the distribution network concerned, while the non-substation buses are connected to loads and renewable energy sources. 
Without loss of generality, assume all the branches $ij \in\mathcal{E}$ are switchable. A binary variable $\alpha_{ij}$ indicates the open ($\alpha_{ij}=0$) and closed ($\alpha_{ij}=1$) status of the switch on branch $ij$.
The uncertain renewable generations and loads are modeled as random variables measured by a joint probability distribution $\boldsymbol{\pi}=(\pi_w,\forall w\in \mathcal{W})$ over a finite set of scenarios $\mathcal{W}$. 

We adapt some widely accepted DNR problem formulations, e.g., those in \cite{peng2014feeder,song2020new}, to the following SDNR problem:
\begin{subequations}\label{eq:SDNR}
\begin{eqnarray}
	&&\textbf{SDNR}:\mathop{\text{min}}_{\boldsymbol{\alpha},\boldsymbol{x}}\sum_{w\in\mathcal{W}}\pi_w\sum_{i\in\mathcal{N}_s}p_{i}^w \label{SDNR_obj}\\
	&&\quad\text{over}\quad \boldsymbol{\alpha}:=(\alpha_{ij},\forall ij\in \mathcal{E})\in \mathcal{A}, \nonumber\\
	&&\quad\quad~~~~~ \boldsymbol{x}:=(\boldsymbol{x}^{w},\forall w\in\mathcal{W}) \nonumber \\
&& \quad\quad\quad~~~~~ = \left(p_{i}^{w}, q_{i}^{w}, \forall i\in \mathcal{N}_s;~v_{i}^{w},\forall i\in\mathcal{N};\right.\nonumber\\
&&\qquad\qquad\qquad \left.  ~p_{ij}^{w},q_{ij}^{w},l_{ij}^{w},\forall ij\in \mathcal{E}; ~\forall w\in\mathcal{W}\right)\nonumber\\
	&&\quad\textnormal{s.t.}\qquad\forall w\in\mathcal{W}:\nonumber\\
	&&\sum_{ij\in\mathcal{E}}\left(p_{ij}^w- \textnormal{r}_{ij}l_{ij}^w\right)+ p_{j}^w=  \sum_{jk\in\mathcal{E}}p_{jk}^w,~\forall j\in\mathcal{N}\label{SDNR_pij}\\
 	&&\sum_{ij\in\mathcal{E}}\left(q_{ij}^w-\textnormal{x}_{ij}l_{ij}^w\right) +q_{j}^w=   \sum_{jk\in\mathcal{E}}q_{jk}^w,~\forall j\in\mathcal{N} \label{SDNR_qij}\\
 	&&v_{i}^{w}-v_{j}^{w}\geq 2(\textnormal{r}_{ij}p_{ij}^{w}+\textnormal{x}_{ij}q_{ij}^{w})-(\textnormal{r}_{ij}^2+\textnormal{x}_{ij}^2)l_{ij}^{w}\nonumber \\
 	&&\qquad\qquad\quad-M(1-\alpha_{ij}),\quad\forall ij\in\mathcal{E}\label{SDNR_v1}\\
 	&&v_{i}^{w}-v_{j}^{w}\leq 2(\textnormal{r}_{ij}p_{ij}^{w}+\textnormal{x}_{ij}q_{ij}^{w})-(\textnormal{r}_{ij}^2+\textnormal{x}_{ij}^2)l_{ij}^{w}\nonumber\\
 	&&\qquad\qquad\quad+M(1-\alpha_{ij}),\quad\forall ij\in\mathcal{E}\label{SDNR_v2}\\
 	&&l_{ij}^{w}v_{i}^{w}= (p_{ij}^{w})^2+(q_{ij}^{w})^2,\quad\forall ij\in\mathcal{E}\label{SDNR_lv}\\
	&&(V_i^{\text{min}})^2\leq v_{i}^{w}\leq(V_i^{\text{max}})^2, \quad\forall i\in\mathcal{N}\label{SDNR_V}\\
	&&p_{i}^{\text{min}}\leq p_{i}^w \leq p_{i}^{\text{max}},~ q_{i}^{\text{min}}\leq q_{i}^w \leq q_{i}^{\text{max}},~\forall i\in\mathcal{N}_{s}\label{SDNR_pqgimax}\\
 	&&(p_{ij}^{w})^2+(q_{ij}^{w})^2\leq (s_{ij}^{\text{max}})^2,\quad \forall ij\in\mathcal{E}\label{SDNR_sijmax}\\
	&& -\alpha_{ij}p_{ij}^{\text{max}}\leq p_{ij}^w \leq \alpha_{ij}p_{ij}^{\text{max}},\quad\forall ij\in\mathcal{E}\label{SDNR_pijmax}\\
	&&  -\alpha_{ij}q_{ij}^{\text{max}}\leq q_{ij}^w \leq \alpha_{ij}q_{ij}^{\text{max}},\quad\forall ij\in\mathcal{E}\label{SDNR_qijmax}\\
 	&&0\leq l_{ij}^{w} \leq \alpha_{ij}(I_{ij}^{\text{max}})^2, \quad\forall ij\in\mathcal{E}. \label{SDNR_lmax}
\end{eqnarray}
\end{subequations}
The notation in \eqref{eq:SDNR} is consistent with that in the classical Dist-Flow model  \cite{baran1989network} and can be found in the Nomenclature, and is thus not repeated here. 
Specifically, the objective \eqref{SDNR_obj} minimizes the expected total active power injection (i.e., power supply from the main grid outside of the distribution network) into all the substation buses. 
Let $\mathcal{G}(\boldsymbol{\alpha})$ denote the network under switch status $\boldsymbol{\alpha}=(\alpha_{ij},\forall ij\in \mathcal{E})$. The feasible set of $\boldsymbol{\alpha}$ is defined as \cite{peng2014feeder}:
\begin{eqnarray}
\mathcal{A}&:=&\left\{\boldsymbol{\alpha}~ |~\mathcal{G}(\boldsymbol{\alpha})~\textnormal{has no loop; and each bus in}~\mathcal{N}_d\right. \nonumber\\
&&\left.\qquad \textnormal{is connected to a single bus in}~\mathcal{N}_s.\right\} \nonumber
\end{eqnarray}
We refer to a network satisfying the condition in $\mathcal{A}$ as a \textit{radial network}.
Equations \eqref{SDNR_pij}\eqref{SDNR_qij} enforce the active and reactive power balance at each bus.  
With big positive constant $M$, inequalities \eqref{SDNR_v1}\eqref{SDNR_v2} become the voltage drop equation across branch $ij$ when $\alpha_{ij}=1$, and otherwise decouple the voltages at buses $i$ and $j$. The quadratic equation \eqref{SDNR_lv} introduces continuous nonconvexity to the SNDR problem, which will be addressed later with relaxation.  
Inequalities \eqref{SDNR_V}--\eqref{SDNR_lmax} impose the operational limits for voltage magnitudes, substation power injections, branch power flows and current magnitudes.

The decision variables in \eqref{eq:SDNR} include power injections $(p_{i}^{w},q_{i}^{w})$ to substation buses $i\in\mathcal{N}_s$ only, while for non-substation buses $i\in\mathcal{N}_d$, the power injections are given as: 
\begin{eqnarray}
&&p_{i}^{w}=\hat{p}_{ri}^{w}-\hat{p}_{di}^{w},\quad q_{i}^{w}=\hat{q}_{ri}^{w}-\hat{q}_{di}^{w},\quad \forall i\in\mathcal{N}_{d} \nonumber
\end{eqnarray}
where $(\hat{p}_{ri}^{w},\hat{q}_{ri}^{w})$ are the active and reactive power generations of the aggregate renewable energy source, and $(\hat{p}_{di}^{w},\hat{q}_{di}^{w})$ are the active and reactive power consumptions of the aggregate load, at bus $i$. As mentioned before, they are uncertain quantities subject to probability distribution $\boldsymbol{\pi}$ over scenarios $w\in\mathcal{W}$.

The SDNR problem \eqref{eq:SDNR} is a typical mixed-integer nonlinear program, which is solvable by an off-the-shelf solver such as Gurobi. However, the numerical results in Section \ref{sec:experiment} for such a solver reveal a substantial room for improvement in terms of optimality and computational efficiency. This is actually what motivated our development of the SDNR method in this paper.

\subsection{SOC Relaxation of OPF}

For each given and fixed switch status $\boldsymbol{\alpha}$, the SDNR problem \eqref{eq:SDNR} is specified as a stochastic optimal power flow (OPF) problem, which is nonconvex due to the quadratic equation \eqref{SDNR_lv}. We adopt the widely used second-order cone (SOC) relaxation \cite{farivar2013branch} to convexify \eqref{SDNR_lv} as:
\begin{eqnarray}
	l_{ij}^{w}v_{i}^{w}&\geq&(p_{ij}^{w})^2+(q_{ij}^{w})^2,~\forall ij\in\mathcal{E} \label{SDNR_lvr}
\end{eqnarray}
which leads to an SOC-relaxed stochastic OPF problem: 
\begin{eqnarray}
	&&\textbf{SOPF-R}(\boldsymbol{\alpha}):~\mathop{\text{min}}\limits_{\boldsymbol{x}}\sum_{w\in\mathcal{W}}\pi_w\sum_{i\in\mathcal{N}_s}p_{i}^w \nonumber\\
	&&\text{s.t. \eqref{SDNR_pij}--\eqref{SDNR_v2}, \eqref{SDNR_V}--\eqref{SDNR_lmax}, \eqref{SDNR_lvr},}~\forall w \in\mathcal{W}\nonumber.
\end{eqnarray}	

Noticing SOPF-R is decoupled across scenarios $w\in\mathcal{W}$, we can solve it as $|\mathcal{W}|$ deterministic problems, each corresponding to one scenario $w\in\mathcal{W}$: 
\begin{eqnarray}
	\textbf{OPF-R}(\boldsymbol{\alpha},w):&&\mathop{\text{min}}\limits_{\boldsymbol{x}^{w}}\sum_{i\in\mathcal{N}_s}p_{i}^{w} \nonumber\\
	&&\text{s.t. \eqref{SDNR_pij}--\eqref{SDNR_v2}, \eqref{SDNR_V}--\eqref{SDNR_lmax}, \eqref{SDNR_lvr}}\nonumber.
\end{eqnarray}	

We assume the technical conditions in \cite{farivar2013branch,low2014convex} are satisfied to make the SOC relaxation exact: the optimal solution of OPF-R attains equality in \eqref{SDNR_lvr}, i.e., it satisfies \eqref{SDNR_lv}. 
The SOPF-R and OPF-R problems introduced here will need to be solved in the proposed method below for SDNR.

\section{Improved Successive Branch Reduction}

\subsection{Improved One-Stage SBR: Open One Redundant Branch}

\begin{figure}[t]
	\centering
	\includegraphics[width=0.9\columnwidth]{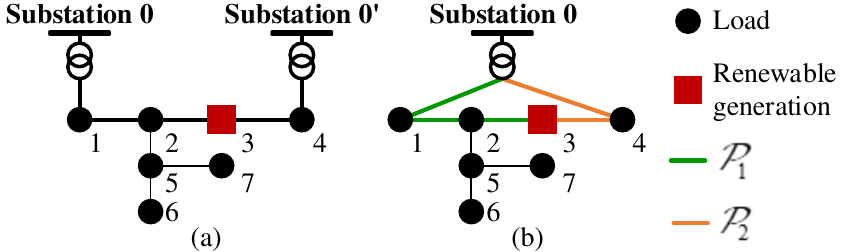}
	\caption{The two substations in (a) can be merged into one substation in (b) without changing the applicability of our method. In (b), buses 0 and 3 inject active power into loop $\mathcal{P}$ and divide $\mathcal{P}$ into two sub-paths $\mathcal{P}_1$ and $\mathcal{P}_2$.}
	\label{fig:DN_one_line_RES}
\end{figure}

We first consider a simple special network that has a single redundant branch compared to a radial network.
According to \cite{peng2014feeder}, this includes two cases as shown in Figure \ref{fig:DN_one_line_RES}: (a) two substations, no loop; and (b) one substation, one loop. Indeed, merging two (or multiple) substations into one does not change the applicability of our method. Therefore, we only consider case (b) when introducing our improved one-stage successive branch reduction (SBR) algorithm below.

Let $\mathcal{P}$ denote the single loop (more precisely the set of branches forming that loop) in the network with a single redundant branch. 
For an arbitrary branch $e = ij \in\mathcal{P}$, let $e^{\text{down}}$ denote its adjacent branch incident to node $j$ (i.e., the downstream branch), and $e^{\text{up}}$ denote its adjacent branch incident to node $i$ (i.e., the upstream branch), both in loop $\mathcal{P}$.  The expectation of the active power flow on branch $ij$ with respect to probability distribution $\boldsymbol{\pi}$ over scenarios $\mathcal{W}$ is denoted as
\begin{eqnarray}
\widetilde{p}_{ij}&:=&\mathbb{E}_{\boldsymbol{\pi}}\left[p_{ij}\right]=\sum_{w\in\mathcal{W}}\pi_{w} p_{ij}^{w}. \nonumber
\end{eqnarray}
Based on this expectation, we define a set of candidate branches that will be used in the improved one-stage SBR algorithm: 
\begin{eqnarray}
	\mathcal{K}(e)&:=&
\begin{cases}
	\{e,e^{\text{down}}\}, &\text{if}~\widetilde{p}_{ij}>0~\text{and}~e^{\text{down}}~\text{exists};\\
	\{e,e^{\text{up}}\}, &\text{if}~\widetilde{p}_{ij}<0~\text{and}~e^{\text{up}}~\text{exists};\\
	\{e\}, &\text{otherwise}.
\end{cases}\label{eq:def:Kappa}
\end{eqnarray}
Further, let $\mathcal{N}_{\mathcal{P}}$ denote the set of buses on loop $\mathcal{P}$. For any bus $i\in\mathcal{N}_{\mathcal{P}}$, define
\begin{eqnarray}
	\widetilde{p}_{i}^{\mathcal{P}} &:=& \mathbb{E}_{\boldsymbol{\pi}}\left[\sum_{\substack{k\notin \mathcal{N}_{\mathcal{P}}\\ki\in\mathcal{E}}}\left(p_{ki}-\textnormal{r}_{ki}l_{ki}\right) +p_{i}-\sum_{\substack{j\notin \mathcal{N}_{\mathcal{P}}\\ij\in\mathcal{E}}}p_{ij}\right]
	\label{eq:pi_equivalent}
\end{eqnarray}
which is the expected active power injection from bus $i$ into loop $\mathcal{P}$. 
Since the total active power loss on all the branches in $\mathcal{P}$ must be supplied by at least one source, there are $n_r\geq 1$ buses $i \in \mathcal{N}_\mathcal{P}$ that satisfy $\widetilde{p}_{i}^{\mathcal{P}}>0$. They divide loop $\mathcal{P}$ into $n_r$ sub-paths $(\mathcal{P}_{m},~m=1,...,n_r)$. Figure \ref{fig:DN_one_line_RES}(b) demonstrates an example where $n_r=2$, and buses 0 and 3 divide loop $\mathcal{P}$ into two sub-paths $\mathcal{P}_1$ and $\mathcal{P}_2$. In particular, active power is injected into loop $\mathcal{P}$ from the main grid (outside of the distribution network) at the substation bus 0, and from the renewable energy source at bus 3.


\begin{algorithm}[t]
	\caption{Improved one-stage SBR}\label{alg:MOSSBR}
	\SetKwData{Left}{left}\SetKwData{This}{this}\SetKwData{Up}{up}\SetKwFunction{Union}{Union}\SetKwFunction{FindCompress}{FindCompress}\SetKwInOut{Input}{Input}\SetKwInOut{Output}{Output}
	Initialize switch status as $\boldsymbol{\alpha}_{\mathcal{E}}:=(\alpha_{ij}=1,\forall ij\in\mathcal{E})$.\\
	Solve $\textbf{OPF-R}(\boldsymbol{\alpha}_{\mathcal{E}},w)$ to obtain the optimal solution $\boldsymbol{x}^w_\mathcal{E}$ for all scenarios $w\in\mathcal{W}$.\\
	Calculate $\widetilde{p}_{i}^{\mathcal{P}}$ by \eqref{eq:pi_equivalent} for each bus $i \in \mathcal{N}_\mathcal{P}$. Buses with $\widetilde{p}_{i}^{\mathcal{P}}>0$ divide loop $\mathcal{P}$ into sub-paths $\{\mathcal{P}_{1},..., \mathcal{P}_{n_r}\}$.\\
	\For{$m=1$ \KwTo $n_r$}
	{
	$\hat{e}_m\leftarrow\text{arg}\mathop{\text{min}}_{e\in\mathcal{P}_{m}}\mathbb{E}_{\boldsymbol{\pi}}\left[|p_{e}(\boldsymbol{x}_{\mathcal{E}})|\right]$
	}
	\For{$e\in \cup_{m=1}^{n_r} \mathcal{K}(\hat{e}_m)$}
	{Solve $\textbf{SOPF-R}(\boldsymbol{\alpha}_{\mathcal{E}\backslash\{e\}})$ and denote its optimal objective value as  $\widetilde{p}_0(e)$.}
	$e^{*}\leftarrow \text{arg}\mathop{\text{min}}_{e\in \cup_{m=1}^{n_r} \mathcal{K}(\hat{e}_m)} \widetilde{p}_0(e)$ \\
	\Return $\boldsymbol{\alpha}^{*}=\boldsymbol{\alpha}_{\mathcal{E}\backslash \{e^{*}\}}$.
\end{algorithm}

Inspired by the baseline SBR algorithm in \cite{peng2014feeder} for deterministic DNR with positive loads only, we propose an improved one-stage SBR algorithm, Algorithm \ref{alg:MOSSBR}, for SDNR with renewable power generations.
Lines 1--2 in Algorithm \ref{alg:MOSSBR} solve the SOC-relaxed stochastic OPF problem $\textbf{SOPF-R}(\boldsymbol{\alpha}_{\mathcal{E}})$ as $|\mathcal{W}|$ deterministic problems $\textbf{OPF-R}(\boldsymbol{\alpha}_{\mathcal{E}},w)$, each corresponding to a scenario $w\in\mathcal{W}$, under switch status $\boldsymbol{\alpha}_{\mathcal{E}}$ that closes all the branches (including the single redundant branch).   
At their optimal solutions, Line 3 calculates the expected active power injections $\widetilde{p}_{i}^{\mathcal{P}}$ for all $i\in\mathcal{N}_\mathcal{P}$, to divide loop $\mathcal{P}$ into sub-paths accordingly.
Using $p_{e}(\boldsymbol{x}_{\mathcal{E}})$ to denote the active power flow on branch $e$ at the optimal solution $\boldsymbol{x}_{\mathcal{E}}=(\boldsymbol{x}_{\mathcal{E}}^w,\forall w\in\mathcal{W})$, Lines 4--6 find the branch $\hat{e}_{m}$ that carries the minimum expected absolute value of active power flow within each sub-path $\mathcal{P}_{m}$. 
The union of $\mathcal{K}(\hat{e}_{m})$ defined in \eqref{eq:def:Kappa} over all the sub-paths $m=1,...,n_r$ is taken as the set of candidate branches to open. For each such candidate branch $e$, Lines 7--9 solve the SOC-relaxed stochastic problem \textbf{SOPF-R} under switch status $\boldsymbol{\alpha}_{\mathcal{E}\backslash\{e\}}$ that opens $e$ only and closes all other branches. Among all such \textbf{SOPF-R} problems solved, the one that attains the minimum optimal objective value determines the optimal branch $e^*$ to be opened by Algorithm \ref{alg:MOSSBR} (Lines 10--11).

\subsection{Rationale behind Algorithm \ref{alg:MOSSBR}}

The rationale behind Algorithm \ref{alg:MOSSBR} is to generalize the baseline SBR algorithm in \cite{peng2014feeder}. Indeed, Algorithm \ref{alg:MOSSBR} is equivalent to the one in \cite{peng2014feeder} in the following idealized special case:
\begin{enumerate}
\item The network has one redundant branch and two substation buses, as shown in Figure \ref{alg:MOSSBR}(a) for instance. Note that we introduced Algorithm \ref{alg:MOSSBR} for a network like Figure \ref{alg:MOSSBR}(b) that has one substation in a single loop $\mathcal{P}$. When the substation is split in two, the loop $\mathcal{P}$ becomes a path $\mathcal{P}$ between the two substations, on which Algorithm \ref{alg:MOSSBR} is still applicable. 

\item A deterministic DNR problem is considered, which can be formulated as a special SDNR problem \eqref{eq:SDNR} with only one scenario $\mathcal{W}=\{w\}$ of probability one. 

\item At the optimal solution of $\textbf{SOPF-R}(\boldsymbol{\alpha}_{\mathcal{E}})$, the active power injection $\widetilde{p}_{i}^{\mathcal{P}}$ into path $\mathcal{P}$ is positive if bus $i$ is one of the two substation buses and negative if bus $i$ is a non-substation bus in $\mathcal{N}_{\mathcal{P}}$. Therefore, Line 3 of Algorithm \ref{alg:MOSSBR} generates only one (trivial) sub-path $\mathcal{P}_1=\mathcal{P}$.   

\item Other technical conditions: constant bus voltage magnitudes, small angle differences across branches, compactness of the feasible set of $\textbf{SOPF-R}(\boldsymbol{\alpha}_{\mathcal{E}})$, and exactness of the SOC relaxation, as elaborated in \cite{peng2014feeder}.
\end{enumerate} 
In this paper, we generalize the special case above in the following two aspects. 

First, in setting 3) above we allow $\widetilde{p}_{i}^{\mathcal{P}}>0$ for non-substation buses $i \in \mathcal{N}_{\mathcal{P}}$ to incorporate active power injections by renewable energy sources. These power-injecting buses divide $\mathcal{P}$ into multiple sub-paths $\{\mathcal{P}_m,~m=1,...,n_r\}$. 
The same procedure as in \cite{peng2014feeder} can prove that the optimal branch to open within each sub-path $\mathcal{P}_m$ must come from the candidate set $\mathcal{K}(\hat{e}_m)$. Therefore, by searching through $\mathcal{K}(\hat{e}_m)$ for all $m=1,...,n_r$ in Lines 7--10, Algorithm \ref{alg:MOSSBR} determines the optimal branch to open for the network. 

Second, we generalize 2) from a single scenario (i.e., a deterministic DNR problem) to multiple (yet a finite number of) scenarios in $\mathcal{W}$. The optimality proof in \cite{peng2014feeder} for the set of variables in the single scenario can be applied in a similar way to their expectations (i.e., finite linear combinations).  

Generalizing 2) and 3) as above and keeping 1) and 4) unchanged, one can prove that Algorithm \ref{alg:MOSSBR} finds the optimal solution to the SDNR problem \eqref{eq:SDNR}. The detailed proof is skipped due to its similarity to the proof in \cite{peng2014feeder}. 
We admit that the technical conditions in 4) are highly idealized and that relaxing them would make the optimality proof much harder. That being said, we hope this discussion may provide some insight and confidence for the application of the proposed algorithm.  

We have to point out that the generalizations above take extra computational efforts associated with the number $n_r$ of active power-injecting buses and the number $|\mathcal{W}|$ of uncertainty scenarios.  
Numerical results in Section \ref{sec:experiment} will show that such increased computations are acceptable given that they produce a satisfactory solution to the difficult SDNR problem with renewable generations.

\subsection{Two-stage SBR: Open Multiple Redundant Branches}

\begin{algorithm}[t]
	\caption{Two-stage SBR}\label{alg:TSSBR}
	\SetKwData{Left}{left}\SetKwData{This}{this}\SetKwData{Up}{up}\SetKwFunction{Union}{Union}\SetKwFunction{FindCompress}{FindCompress}\SetKwInOut{Input}{Input}\SetKwInOut{Output}{Output}
	\textbf{First stage:}\\
	Initialize switch status as $\boldsymbol{\alpha}_{\mathcal{E}}:=(\alpha_{ij}=1,\forall ij\in\mathcal{E})$.\\
	Solve $\textbf{OPF-R}(\boldsymbol{\alpha}_{\mathcal{E}},w)$ to obtain the optimal solution $\boldsymbol{x}^w_\mathcal{E}$ for all scenarios $w\in\mathcal{W}$.\\
	\For{$l=1$ \KwTo $L$}
	{
	$e_{l}^{o}\leftarrow\text{arg}\mathop{\text{min}}_{e\in\mathcal{P}^{l}}\mathbb{E}_{\boldsymbol{\pi}}\left[|p_{e}(\boldsymbol{x}_{\mathcal{E}})|\right]$\\
	Open branch $e_{l}^{o}$.\\
	\If{$e_{l}^{o}$ is a common branch of $\mathcal{P}^{l}$ and $\mathcal{P}^{k}$}{
       Update loop $\mathcal{P}^{k}$ with $e_{l}^{o}$ open.
      }	
	}
    $\mathcal{E}_o\leftarrow\{e_{l}^{o},~l=1,...,L\}$\\
	\textbf{Second stage:}\\
	\For{$l=1$ \KwTo $L$}
	{
	Define a set of open branches $\mathcal{E}_o^l:=\mathcal{E}_o\backslash\{e_l^o\}$.\\
	Run Algorithm \ref{alg:MOSSBR} with initial switch status $\boldsymbol{\alpha}_{\mathcal{E}\backslash\mathcal{E}_o^l}$, to obtain an optimal branch $e_l^*$ to open; denote the corresponding optimal objective value as $\widetilde{p}_0(e_l^*)$.\\
	}
	$l_{\text{min}}\leftarrow \text{arg}\mathop{\text{min}}_{l=1,...,L} \widetilde{p}_0(e_l^*)$ \\
	\Return $\boldsymbol{\alpha}^{*}=\boldsymbol{\alpha}_{\mathcal{E}\backslash(\mathcal{E}_o^{l_{\text{min}}}\cup\{e_{l_{\text{min}}}^*\})}$
\end{algorithm}

\begin{figure}[t]
	\centering
	\includegraphics[width=0.9\columnwidth]{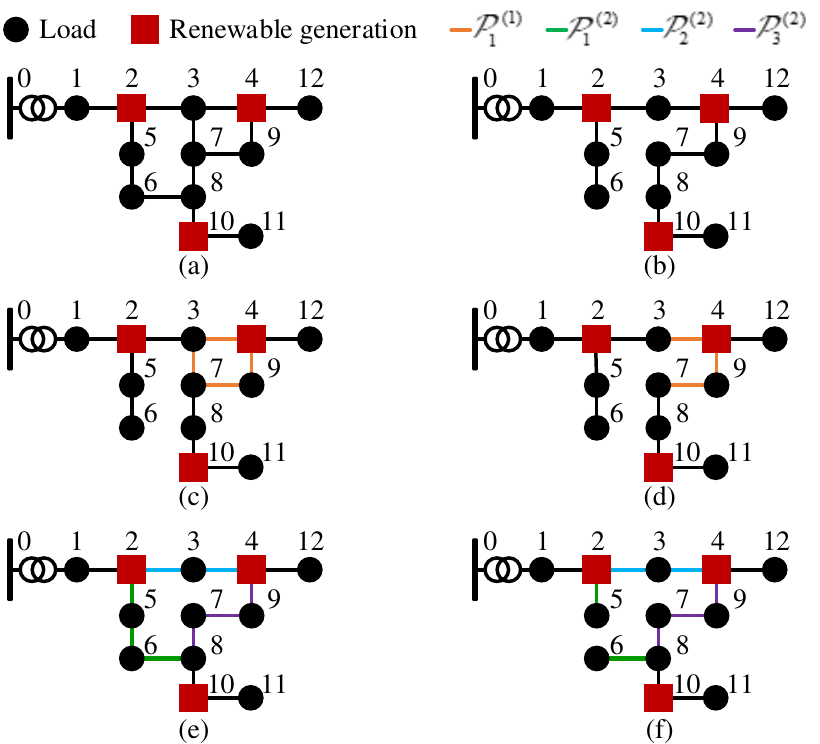}
	\caption{Illustration of the two-stage SBR (Algorithm \ref{alg:TSSBR}). First stage: (a) the initial network with two loops; (b) open branches $\mathcal{E}_o=\{e_1^o, e_2^o\}=\{(3,7), (6,8)\}$. Second stage: (c) close branch $e_1^o=(3,7)$ to form loop $\mathcal{P}^{(1)}$, which is divided into only one sub-path $\mathcal{P}^{(1)}_1$ by active power-injecting bus 4; (d) run Algorithm \ref{alg:MOSSBR} to open branch $e_1^{*}=(3,7)$; (e) close branch $e_2^o=(6,8)$ to form loop $\mathcal{P}^{(2)}$, which is divided into three sub-paths $\mathcal{P}^{(2)}_m$, $m=1,2,3$ by active power-injecting buses 2, 4, 8; (f) run Algorithm \ref{alg:MOSSBR} to open branch $e_2^{*}=(5,6)$. Radial networks (d) and (f) need to be compared, and the one that attains the lower optimal objective value in Algorithm \ref{alg:MOSSBR} is the final output of Algorithm \ref{alg:TSSBR}.}
	\label{fig:TSSBR_eg}
\end{figure}

We now consider a general network with $L>1$ redundant branches that form $L$ chordless loops $\mathcal{P}^{l}$, $l=1,...,L$. 
Without loss of generality, we still assume the network has a single substation bus (that may merge multiple substations). To solve the SDNR problem in such a network, we propose a heuristic two-stage SBR algorithm, Algorithm \ref{alg:TSSBR}, which iteratively calls Algorithm \ref{alg:MOSSBR} as will be elaborated below.

The first stage of Algorithm \ref{alg:TSSBR} starts from switch status $\boldsymbol{\alpha}_{\mathcal{E}}$ that closes all the branches, including $L$ redundant branches (Line 2). 
Line 3 solves the SOC-relaxed stochastic OPF problem $\textbf{SOPF-R}(\boldsymbol{\alpha}_{\mathcal{E}})$ as $|\mathcal{W}|$ deterministic problems $\textbf{OPF-R}(\boldsymbol{\alpha}_{\mathcal{E}},w)$, each corresponding to a scenario $w\in\mathcal{W}$, under the initial switch status $\boldsymbol{\alpha}_{\mathcal{E}}$.
Lines 4--5 find the branch $e_l^o$ that carries the minimum expected absolute value of active power flow $p_{e}(\boldsymbol{x}_{\mathcal{E}})$ at the optimal solution $\boldsymbol{x}_{\mathcal{E}}$, within each loop $\mathcal{P}^{l}$. Line 6 opens that branch $e_l^o$. If branch $e_l^o$ simultaneously lies in two loops including $\mathcal{P}^{l}$, opening $e_l^o$ will also change the other loop, say $\mathcal{P}^{k}$ (Lines 7--9). The open branches $\{e_l^o,~l=1,...,L\}$ are collected in a 
set $\mathcal{E}_o$ (Line 11), which serves as the basis for the second stage of Algorithm \ref{alg:TSSBR}.

For illustration, Figure \ref{fig:TSSBR_eg}(a) shows an example with $L=2$, where loop $\mathcal{P}^1=\{(2,3), (3,7), (7,8), (2,5), (5,6), (6,8)\}$ and $\mathcal{P}^2 = \{(3,4), (4,9), (3,7), (7,9)\}$. Branch $(3,7)$ is a common branch of the two loops. 
In the first stage of Algorithm \ref{alg:TSSBR}, branch $e_1^o=(3,7)$ is first opened in loop $\mathcal{P}^{1}$, which causes loop $\mathcal{P}^{2}$ to change into a larger loop; then branch $e_2^o=(6,8)$ is opened in the updated loop $\mathcal{P}^{2}$, leading to the radial network in Figure \ref{fig:TSSBR_eg}(b). The set of open branches obtained in the first stage is $\mathcal{E}_o=\{e_1^o, e_2^o\}=\{(3,7), (6,8)\}$.  

The second stage of Algorithm \ref{alg:TSSBR} runs $L$ iterations. In each iteration $l=1,...,L$, all the branches in $\mathcal{E}_o$ (obtained in the first stage) except $e_l^o$ are open, while all other branches including $e_l^o$ are closed (Lines 13--14). This switch status, denoted as $\boldsymbol{\alpha}_{\mathcal{E}\backslash\mathcal{E}_o^l}$, has only one redundant line, and thus Algorithm \ref{alg:MOSSBR} can be applied to it to obtain an optimal switch $e_l^*$ to open (Line 15). Note that $e_l^*$ may or may not be the same branch as $e_l^o$. In this way, each of the $L$ iterations leads to a radial network. Among the $L$ radial networks obtained, the one whose optimal objective value in Algorithm \ref{alg:MOSSBR} attains the minimum, denoted as the $l_{\text{min}}$-th one, decides the ultimate switch status to be returned by Algorithm \ref{alg:TSSBR} (Lines 17--18).   


In Figure \ref{alg:TSSBR}, the second stage of Algorithm \ref{alg:TSSBR} is run for two iterations $l=1,2$. 
In the first iteration, branch $e_1^o=(3,7)$ is closed to form 
a single loop $\mathcal{P}^{(1)}$ in Figure \ref{alg:TSSBR}(c); then Algorithm \ref{alg:MOSSBR} returns an optimal branch $e_1^*=(3,7)$ to open, as Figure \ref{fig:TSSBR_eg}(d) shows. In the second iteration, branch $e_2^o=(6,8)$ is closed to form 
a single loop $\mathcal{P}^{(2)}$ in Figure \ref{alg:TSSBR}(e); then Algorithm \ref{alg:MOSSBR} returns an optimal branch $e_2^*=(5,6)$ to open, as Figure \ref{fig:TSSBR_eg}(f) shows. Finally, the optimal objective values attained by Algorithm \ref{alg:MOSSBR} need to be compared between \ref{fig:TSSBR_eg}(d) and \ref{alg:TSSBR}(f), and the lower one decides the optimal switch status $\boldsymbol{\alpha}^*$.

Both Algorithm \ref{alg:TSSBR} and the heuristic \cite[Algorithm 3]{peng2014feeder} are designed for a general network with multiple redundant branches. Besides the stochastic setting and the active renewable power injections handled in Algorithm \ref{alg:TSSBR} as discussed before, their major difference lies in the iterative close-and-open procedure in Lines 13--16 of Algorithm \ref{alg:TSSBR}, which further improves the optimality of SDNR solution compared to \cite[Algorithm 3]{peng2014feeder} as shown numerically in Section \ref{sec:experiment}. 
We propose Algorithm \ref{alg:TSSBR} as an experimentally effective and research-wise insightful heuristic. Its formal performance proof appears to be much harder and is left beyond this work.  


\section{Case Studies}\label{sec:experiment}

\subsection{Experimental Setup}

\begin{figure}[t]
	\centering
	\includegraphics[width=0.95\columnwidth]{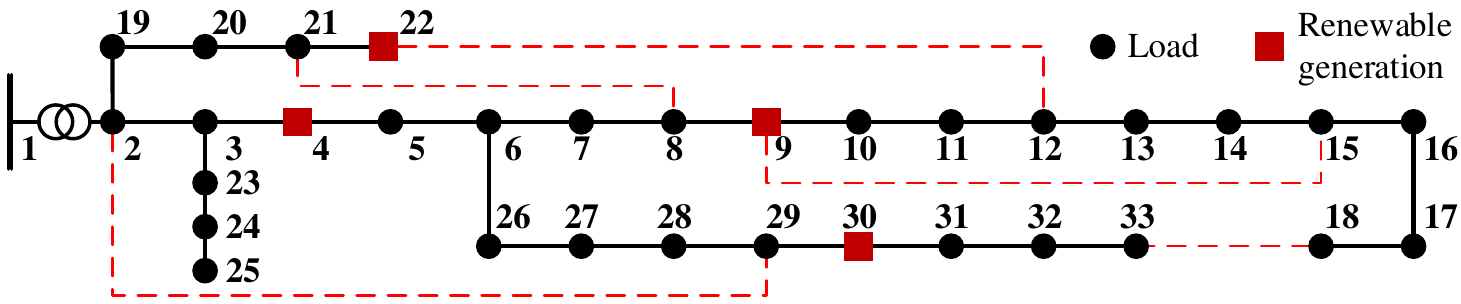}
	\caption{IEEE 33-bus distribution network with five initial redundant branches (red dashed lines).}
	\label{fig:IEEE_33}
\end{figure} 

\begin{figure}[t]
	\centering
	\includegraphics[width=0.99\columnwidth]{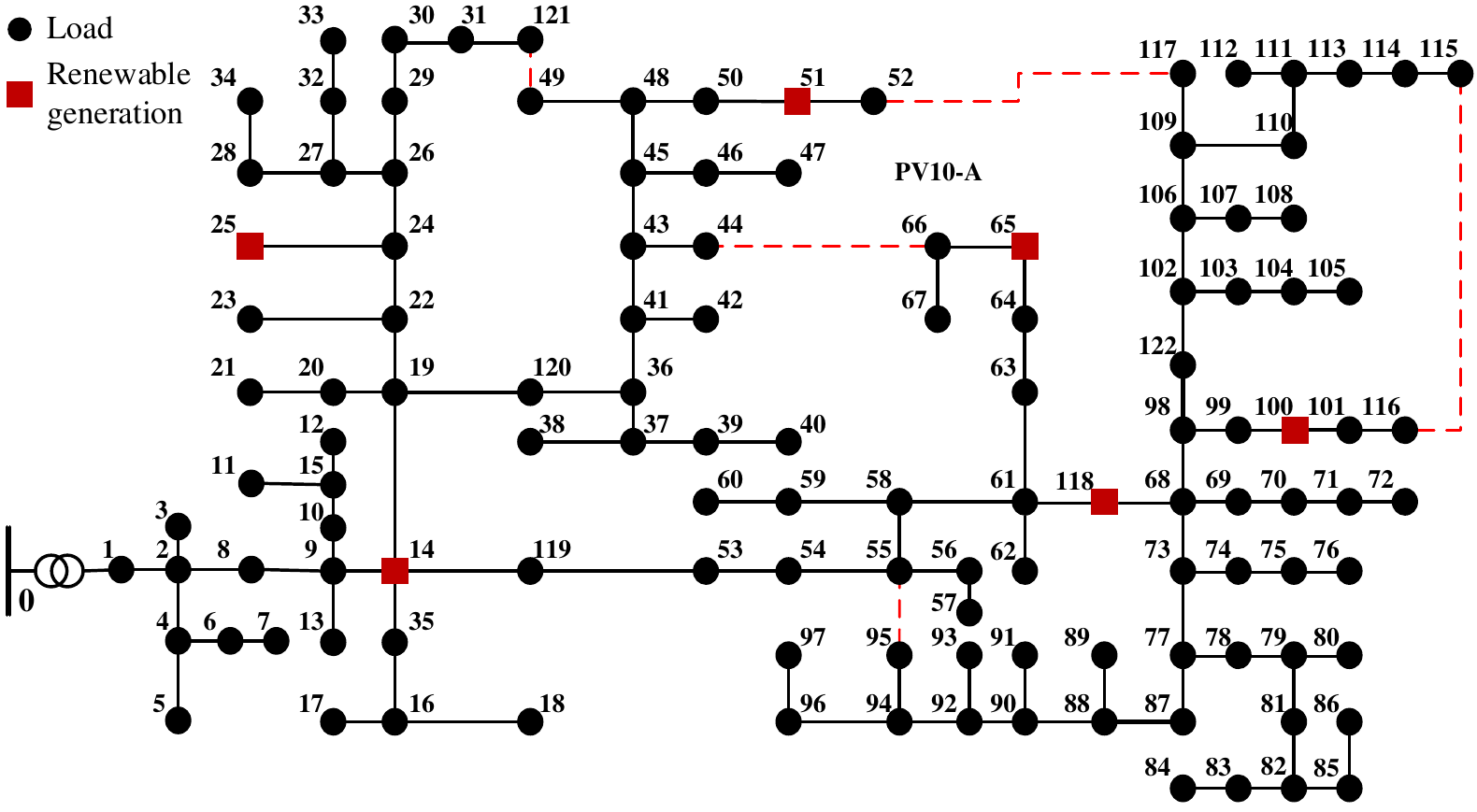}
	\caption{IEEE 123-bus distribution network with five initial redundant branches (red dashed lines).}
	\label{fig:IEEE_123}
\end{figure} 

To validate the proposed method, we conduct numerical experiments with the IEEE 33-bus and 123-bus distribution network models. 
As shown in Figures \ref{fig:IEEE_33} and \ref{fig:IEEE_123}, each network has one substation bus and five redundant branches. 
Each non-substation bus connects to either a load or a pair of renewable energy sources (a small wind turbine and a solar panel). 
We use the load, wind and solar generation data in Germany from April to June 2020 in hourly resolution \cite{OPSDP2022}, scale them to fit the IEEE network capacities, and cluster them into a certain number of scenarios (depending on the test case) using $k$-medoids \cite{huang2021distribution}. 
After getting the active power of loads and renewable generations, we determine their reactive power by fixed power factors \cite{zhan2020switch}. 

 The experiments are run on a 64-bit MacBook with Apple M1 Pro Chip, 8-core CPU, and 32GB RAM. In the proposed method, we can solve the SOC-relaxed problems SOPF-R and OPF-R using a convex solver as currently written in Algorithms \ref{alg:MOSSBR} and \ref{alg:TSSBR}, or solve the OPF problems without relaxation using MATPOWER. In these experiments, we choose the latter.  

The proposed method is compared with the following two methods in terms of optimality and computational efficiency. 

Method 1: Mathematical programming. We use the commercial mixed-integer nonlinear solver Gurobi  
	to solve SDNR in its original form \eqref{eq:SDNR} except that \eqref{SDNR_lv} is replaced by its SOC relaxation \eqref{SDNR_lvr} to reduce computational burden to some extent.
	
Method 2: Stochastic versions of the baseline algorithms, i.e., \cite[Algorithm 1]{peng2014feeder} to open a single redundant branch and \cite[Algorithm 3]{peng2014feeder} to open multiple redundant branches. 
	The deterministic active power flows in the baseline algorithms are replaced by their expectations. Particularly, compared to the proposed method, \cite[Algorithm 1]{peng2014feeder} did not divide the loop into sub-paths at the active power-injecting buses, and \cite[Algorithm 3]{peng2014feeder} did not apply the close-and-open procedure in Lines 13--16, Algorithm \ref{alg:TSSBR}. 
	
	Other methods, such as the switch opening and exchange \cite{zhan2020switch} and genetic algorithm \cite{jakus2020optimal}, are not compared in our experiments because of their obviously more complicated implementations and heavier computational burdens than the proposed method. 


\subsection{Optimality of Algorithm \ref{alg:MOSSBR}: One-Stage SBR}

\begin{table}[t]
	\centering
	\caption{Relative errors in minimum expected total power losses (IEEE 33-bus network with a single redundant branch)}
	\renewcommand\arraystretch{1.2}	
	\begin{tabular}{c|ccc|ccc}		
		\hline
		\hline 
		\multirow{3}{*}{\makecell{Initial\\redundant branch}}  
		&\multicolumn{3}{c|}{\makecell{Algorithm \ref{alg:MOSSBR} (proposed)\\vs. Method 1}}  &\multicolumn{3}{c}{\makecell{Method 2\\vs. Method 1}}\\
		\cline{2-7}
		&mean  &max &min  &mean  &max &min \\	
		&(\%)  &(\%) &(\%)  &(\%)  &(\%) &(\%)\\	
		\hline 
		$(8,21)$ &0.04 &0.29 &-0.05 &0.04  &0.29 &-0.05\\
		\hline
		$(9,15)$ &-0.02 &0.16 &-0.06 &-0.02 &0.16 &-0.06\\
		\hline
		$(12,22)$ &-0.11 &0 &-0.82 &0.30  &0.83 &-0.31\\
		\hline
		$(18,33)$ &-0.21 &0 &-1.17 &-0.21  &0 &-1.17\\
		\hline
		$(2,29)$  &0 &0 &0 &0 &0 &0\\
		\hline  
		\hline  	
	\end{tabular} \label{tab_MOSSBR}
\end{table}

\begin{figure}[t]
	\centering
	\includegraphics[width=0.95\columnwidth]{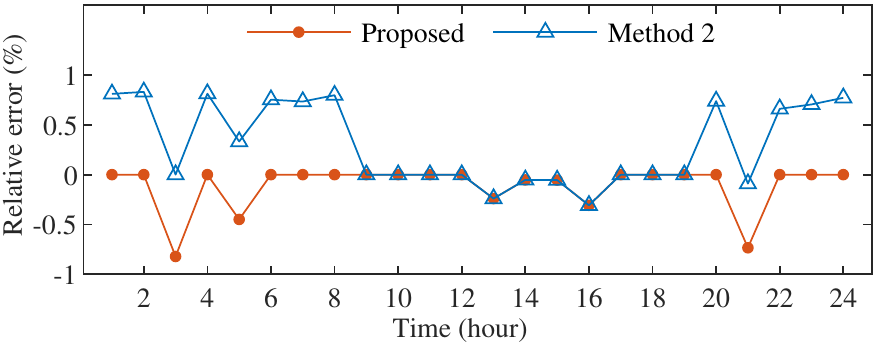}
	\caption{Relative errors in minimum expected total power losses, in the IEEE 33-bus network with initial redundant branch $(12,22)$. The proposed algorithm refers to Algorithm \ref{alg:MOSSBR}, the improved one-stage SBR.}
	\label{fig:MOSSBR_24hrs}
\end{figure}

We first experiment with the improved one-stage SBR algorithm, Algorithm \ref{alg:MOSSBR}, in the IEEE 33-bus network. 
Each of the five redundant branches in Figure \ref{fig:IEEE_33} is closed to create a network with a single (initial) redundant branch.  
On each of the five networks created, we show in Table \ref{tab_MOSSBR} the relative errors in minimum expected total power losses between Algorithm \ref{alg:MOSSBR}  (the proposed) and Method 1 in columns 1--3, and between Method 2 and Method 1 in columns 4--6. 
Here Method 1, i.e., mathematical programming with Gurobi, is used as a trustworthy benchmark to define the relative errors for other methods. 
The mean, maximum, and minimum records over the 24 hourly SDNR problems in a day are listed in the table.

It is observed from Table \ref{tab_MOSSBR} that the relative errors in minimum expected total power losses of Algorithm \ref{alg:MOSSBR} are mostly negligible, with the largest error just around $1\%$.
Actually, Algorithm \ref{alg:MOSSBR} and Method 2 find the same SDNR solutions in each of the five networks except the one with initial redundant branch $(12,22)$; for $(12,22)$, Algorithm \ref{alg:MOSSBR} finds the same or better solutions than Method 2 over 24 hours, as shown in Figure \ref{fig:MOSSBR_24hrs}. In general, the optimality of Algorithm \ref{alg:MOSSBR} is as competitive as its peers.


\subsection{Optimality of Algorithm \ref{alg:TSSBR}: Two-Stage SBR}

\begin{table}[t]
	\centering
	\caption{Relative errors in minimum expected total power losses under different renewable penetration levels (five redundant branches)}
	\renewcommand\arraystretch{1.2}	
	\begin{tabular}{c|c|ccc|ccc}	
		\hline
		\hline 
		\multirow{3}{*}{$k_{r}$} 
		&\multirow{3}{*}{\makecell{Method}}  &\multicolumn{3}{c|}{33-bus network}  &\multicolumn{3}{c}{123-bus network}\\
		\cline{3-8}
		& &mean  &max &min  &mean  &max &min \\	
		& &(\%)  &(\%) &(\%)  &(\%)  &(\%) &(\%)\\	
		 \hline 
		\multirow{3}{*}{0.5}
		&Proposed &-1.22 &0 &-5.89 &-4.25 &-0.92 &-8.62\\
		&Method 2 &-0.49 &1.28 &-4.70 &-4.04 &-0.70 &-8.42\\
		\hline 
		\multirow{3}{*}{1.0}
		&Proposed &-1.17 &0.21 &-4.06 &-4.52 &-1.43 &-8.59\\
		&Method 2 &0.51 &6.88 &-3.60 &-3.40 &0.26 &-8.39\\
		\hline
		\multirow{3}{*}{2.0}
		&Proposed &-1.80 &0.93 &-4.86 &-5.11 &-0.51 &-9.40\\
		&Method 2 &16.78 &69.49 &-2.83 &-4.16 &1.90 &-7.67\\
		\hline  	
		\multirow{3}{*}{3.0}
		&Proposed &-0.84 &6.08 &-4.58 &-3.89 &0.18 &-8.99\\
		&Method 2 &12.27 &23.84 &1.04 &5.50 &23.42 &-8.80\\	
		\hline  
		\hline  	
	\end{tabular} \label{tab_TSSBR_RES}
\end{table}

\begin{figure}[t]
	\centering
	\includegraphics[width=0.95\columnwidth]{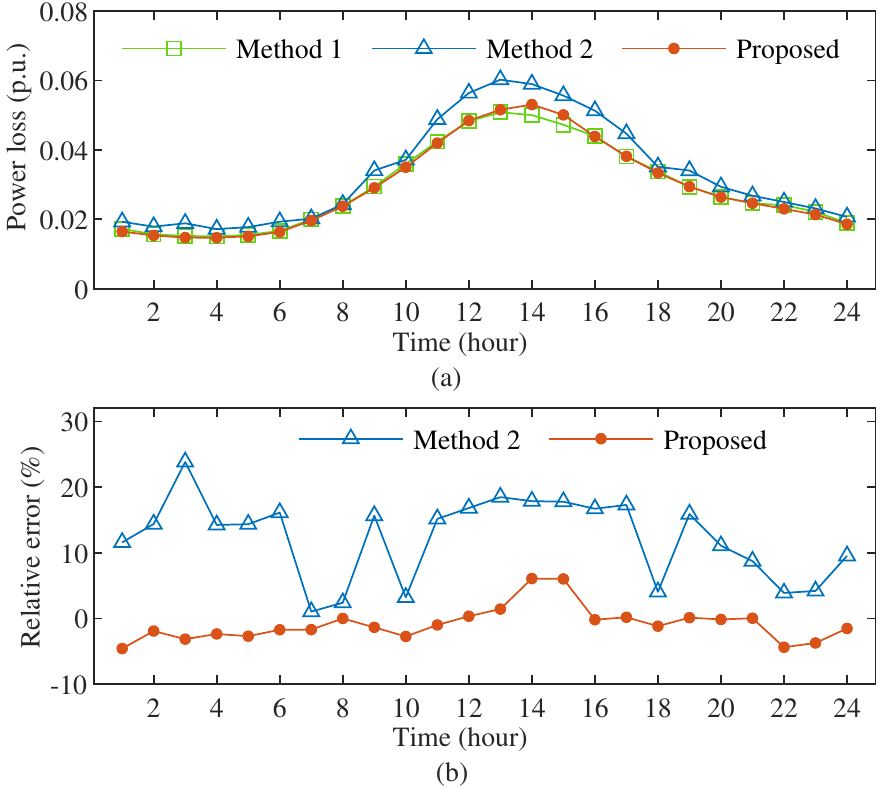}
	\caption{(a) The minimum expected total power losses and (b) the relative errors in minimum expected total power losses, obtained by different methods in the IEEE 33-bus network. Renewable penetration $k_{r}=3$. The proposed algorithm refers to Algorithm \ref{alg:TSSBR}, the two-stage SBR.}
	\label{fig:TSSBR_case33_kres3}
\end{figure}

\begin{figure}[t]
	\centering
	\includegraphics[width=0.95\columnwidth]{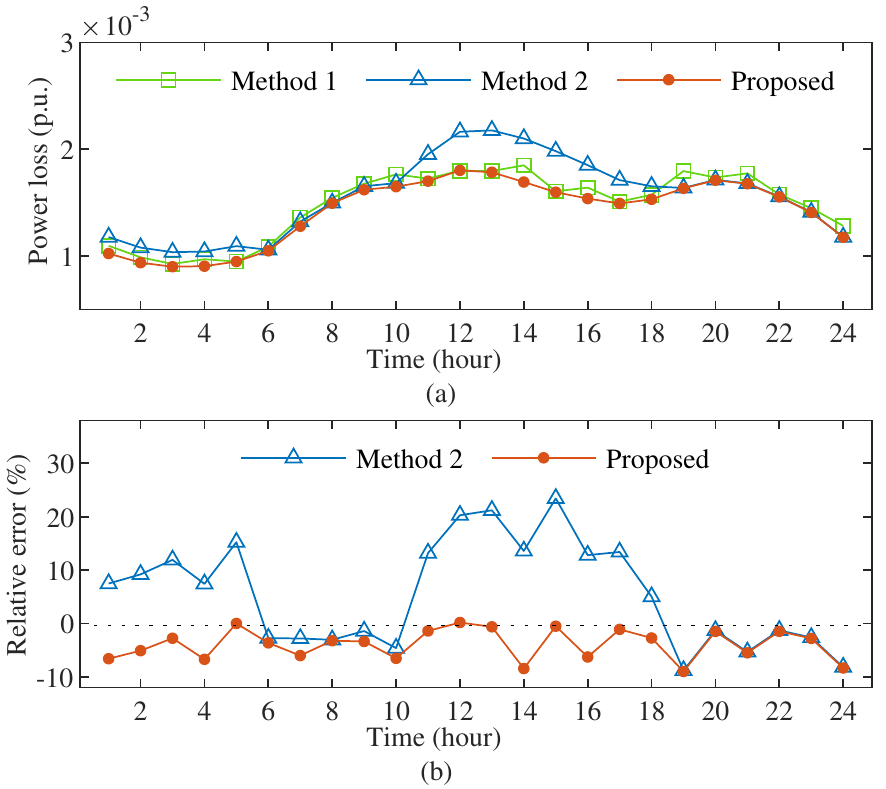}
	\caption{(a) The minimum expected total power losses and (b) the relative errors in minimum expected total power losses, obtained by different methods in the IEEE 123-bus network. Renewable penetration $k_{r}=3$. The proposed algorithm refers to Algorithm \ref{alg:TSSBR}, the two-stage SBR.}
	\label{fig:TSSBR_case123_kres3}
\end{figure}

We now show the performance of the proposed two-stage SBR algorithm, Algorithm \ref{alg:TSSBR}, in the IEEE 33-bus and 123-bus networks, each having five redundant branches. 

The first group of results, collected in Table \ref{tab_TSSBR_RES}, are the relative errors (compared to Method 1) in minimum expected total power losses, under different renewable penetration levels. The renewable penetration level is represented by coefficient $k_{r}$, which scales the renewable generation capacities in the IEEE 33-bus and 123-bus networks. 
An immediate observation is that the proposed algorithm (Algorithm \ref{alg:TSSBR}) attains lower (better) objective values in the larger 123-bus network than the 33-bus network. 
Furthermore, we see that the proposed algorithm attains better objective values than Methods 1 and 2 under low renewable penetration levels $k_{r}$. As $k_{r}$ increases, the optimality of Method 2 severely degrades, while the proposed algorithm still finds better solutions than Method 1 in most cases. For a clearer comparison, we display (a) the minimum expected total power losses and (b) the relative errors in minimum expected total power losses, over 24 hours under a high renewable penetration level $k_{r}=3$, in Figures \ref{fig:TSSBR_case33_kres3} and \ref{fig:TSSBR_case123_kres3} for the 33-bus and 123-bus networks, respectively. These figures validate our conclusion from Table \ref{tab_TSSBR_RES} that the proposed algorithm outperforms Method 2 in optimality under high renewable penetration levels. 


\begin{table}[t]
	\centering
	\caption{Relative errors in minimum expected total power losses under different numbers of scenarios (33-bus with five redundant branches)}
	\renewcommand\arraystretch{1.2}	
	\begin{tabular}{c|c|ccc}		
		\hline
		\hline 
		$|\mathcal{W}|$
		&Method &mean (\%)  &max (\%) &min (\%)  \\		
		 \hline 
		\multirow{3}{*}{5}
		&Proposed &-1.17 &0.21 &-4.06 \\
		\cline{2-5}
		&Method 2 &0.51 &6.88 &-3.60 \\
		\hline 
		\multirow{3}{*}{20}
		&Proposed &-2.93 &-0.50 &-8.13 \\
		\cline{2-5}
		&Method 2 &-2.59 &-0.14 &-7.71 \\
		\hline
		\multirow{3}{*}{40}
		&Proposed &-11.03 &-1.77 &-30.19 \\
		\cline{2-5}
		&Method 2 &-10.71 &-1.77 &-29.83 \\
		\hline  
		\hline  	
	\end{tabular} \label{tab_TSSBR_scenario}
\end{table}

\begin{figure}[t]
	\centering
	\includegraphics[width=0.95\columnwidth]{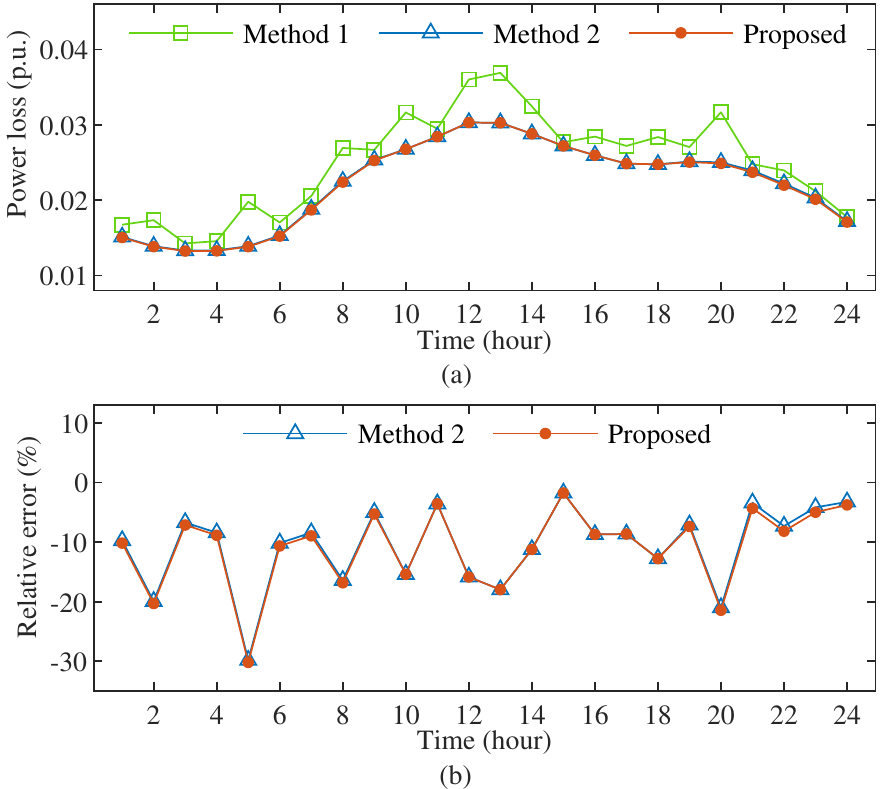}
	\caption{(a) The minimum expected total power losses and (b) the relative errors in minimum expected total power losses, obtained by different methods in the IEEE 33-bus network. The number of uncertainty scenarios $|\mathcal{W}|=40$. The proposed algorithm refers to Algorithm \ref{alg:TSSBR}, the two-stage SBR.}
	\label{fig:TSSBR_case33_Ns40}
\end{figure}

The second group of results, collected in Table~\ref{tab_TSSBR_scenario}, are the relative errors (compared to Method 1) in minimum expected total power losses, under different numbers $|\mathcal{W}|$ of the scenarios of uncertain renewable generations and loads. 
As $|\mathcal{W}|$ increases, both Algorithm 2 (the proposed) and Method 2 find higher-quality solutions. This is also verified by Figure \ref{fig:TSSBR_case33_Ns40} showing (a) the minimum expected total power losses and (b) the relative errors in minimum expected total power losses, over 24 hours under the number of uncertainty scenarios $|\mathcal{W}|=40$.  

A general observation so far is that the proposed SBR method finds better solutions as the network size, the number of redundant lines, and the number of uncertainty scenarios increase. 
Our conjecture is that the increasing numbers of constraints and decision variables (especially the binary switch variables) make the mixed-integer nonlinear SDNR problem more challenging to solve, thus making the superiority of the proposed method more obvious in terms of optimality.  

\subsection{Computational Efficiency}

\begin{figure}[t]
	\centering
	\includegraphics[width=0.95\columnwidth]{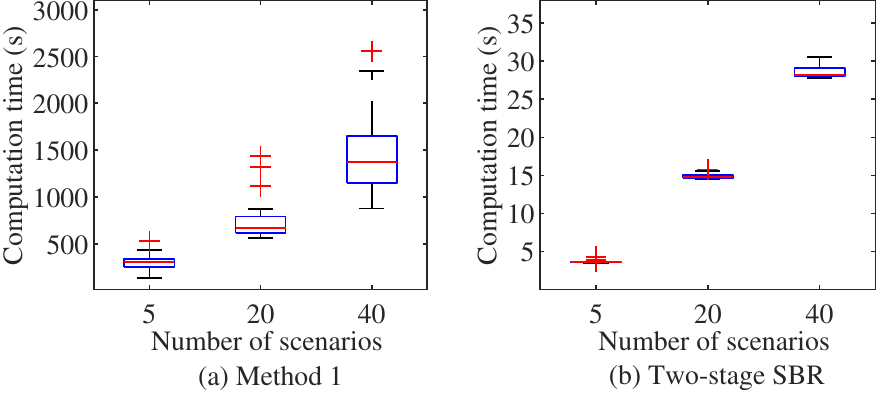}
	\caption{Computation time for SDNR using two methods, in the IEEE 33-bus network under different numbers of uncertainty scenarios.}
	\label{fig:case33_time}
\end{figure}

\begin{figure}[t]
	\centering
	\includegraphics[width=0.95\columnwidth]{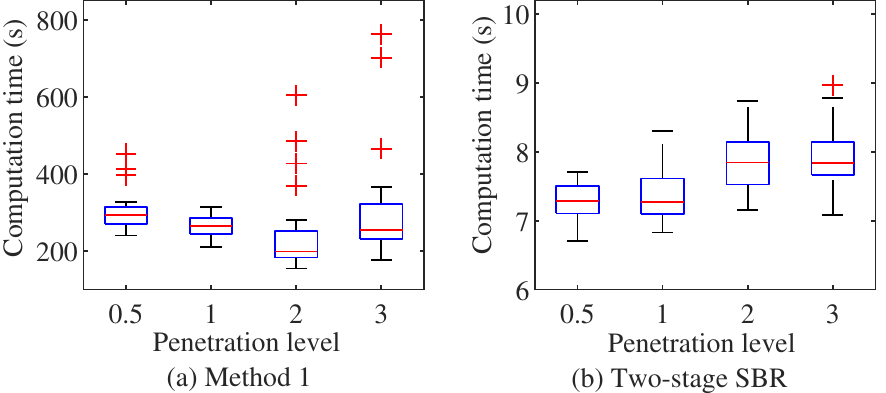}
	\caption{Computation time for SDNR using two methods, in the IEEE 123-bus network under different renewable penetration levels.}
	\label{fig:case123_time}
\end{figure}

Figure \ref{fig:case33_time} plots the computation time (the median, the 25th and 75th percentiles, and a few outliers from the 24 samples in a day) for Method 1 and the proposed two-stage SBR algorithm (Algorithm \ref{alg:TSSBR}) in the IEEE 33-bus network under different numbers of scenarios $|\mathcal{W}|$. Noticing the drastic difference between the vertical axes of the two methods, we see the proposed algorithm is computationally more efficient than Method 1 that uses Gurobi. 

The advantage of the proposed algorithm in computational efficiency is also shown in Figure \ref{fig:case123_time}, which compares its computation time with Method 1 in the IEEE 123-bus network under varying renewable penetration levels $k_r$.

\section{Conclusion}
This paper proposed an improved successive branch reduction (SBR) method to solve stochastic distribution network reconfiguration (SDNR) considering uncertain renewable generations and loads. 
First, for a simple special network with a single redundant branch, we developed an improved one-stage SBR algorithm to incorporate uncertain renewable power generations. 
Then, for a general network with multiple redundant branches, we proposed a two-stage SBR algorithm featuring an iterative close-and-open procedure that runs the one-stage SBR algorithm in each iteration. 
Numerical experiments on the IEEE 33-bus and 123-bus network models validated the improved optimality and computational efficiency of the proposed method compared to a common mixed-integer nonlinear solver and a baseline SBR method from the literature.

In the future, we plan to formally analyze the optimality of the proposed two-stage SBR algorithm (which is currently a heuristic), hoping to provide some insight for an improved SDNR algorithm design with performance guarantees. Another direction of our ongoing work is to develop an efficient robust DNR algorithm incorporating other important considerations such as the small-signal or short-term voltage stability.  

\bibliographystyle{IEEEtran}
\bibliography{references}

\end{document}